\documentstyle[11pt]{article}
\newcommand{\proof}{\noindent {\bf Proof: }}

\newcommand{\corollary}{\noindent {\bf Corollary:}}

\newtheorem{theorem}{Theorem}
\newtheorem{lemma}{Lemma}
\newtheorem{statement}{Statement}
\newtheorem{defi}{Definition}
\def\qed{\hfill $\Box$}

\begin{document}
\title{On the shadow boundary of a centrally symmetric convex body}
\author{\'A.G.Horv\'ath\\ Department of Geometry, \\
Technical University of Budapest,\\
H-1521 Budapest,\\
Hungary}
\date{Febr. 20, 2007}

\maketitle

\section{Introduction, historical remarks}

If $K$ is a $0$-symmetric, bounded, convex body in the Euclidean
$n$-space $R^n$ (with a fixed origin O) then it defines a norm whose
unit ball is $K$ itself (see \cite{l-g}). Such a
space is called {\bf Minkowski normed space}. The main results in
this topic collected in the survey \cite{martini-swanepoel 1} and
\cite{martini-swanepoel 2}. In fact, the norm is a continuous
function which is considered (in the geometric terminology as in
\cite{l-g}) as a gauge function. The metric (the so-called Minkowski
metric), the distance of two points, induced by this norm, is
invariant with respect to the translations of the space.

The unit ball is said to be {\bf strictly convex} if its boundary
contains no line segment.

In previous papers of this topic (\cite{gho1}, \cite{gho2}), we
examined the boundary of the unit ball of the norm and gave two
theorems (Theorem 2 and Theorem 3) similar to the characterization
of the Euclidean norm investigated by H.Mann, A.C.Woods and
P.M.Gruber in \cite{mann}, \cite{woods}, \cite{gruber 1},
\cite{gruber 2} and \cite{gruber 3}, respectively. We proved that if
the unit ball of a Minkowski normed space is strictly convex then
every bisector (which is the collection of those points of the
embedding Euclidean space which have the same distance with respect
to the Minkowskian norm to two given points of the space) is a
topological hyperplane (Theorem 2). Example 3 in \cite {gho1} showed
that strict convexity does not follows from the fact that all
bisectors are topological hyperplanes.

We examined the connection between the shadow boundaries of the unit
ball and the bisectors of the Minkowskian space. We were sure that
the following statement is true: {\bf A bisector is a  topological
hyperplane (meaning that a homeomorphism of $R^{n}$ onto itself
sends the bisector onto a standard hyperplane of $R^{(n-1)}$ ) if
and only if the corresponding shadow boundary is a topological
$(n-2)$-dimensional sphere}, however we proved the conjecture only in
the three dimensional case. (Theorem 2 and Theorem 4) We examined
the basic properties of the shadow boundary (Section 2) and define a
well-usable class of sets - the so-called general parameter spheres.

In this paper we discuss some further topological observations on
shadow boundary and general parameter spheres, we prove that in
general they are not ANR (Absolute neighbourhood retract, see in
\cite{borsuk} or \cite{Repovs} ), but are compact metric spaces,
containing an $(n-2)$-dimensional closed, connected subset separating
the boundary of $K$. Using the approximation theorem of cell-like
mappings we also investigate the manifold case proving that in this
case they are homeomorphic to the $(n-2)$-dimensional sphere. A
consequence of this result that if the bisector is a homeomorphic
copy of $R^{(n-1)}$ then the shadow boundary is a topological
$(n-2)$-sphere, proving the first
direction of the conjecture above. We have two question
corresponding to this conjecture. The converse statement is true or
not?
Is it possible that in the manifold case the embedding of the
bisector and the shadow boundary are not standard ones? In the fourth
paragraph we prove that the embedding of the examined sets (in the
manifold case) are always standard, however the first question is open.

\section{Once more again on the shadow boundary of the
unit ball.}

There are several known properties of the shadow boundary of a
convex body with respect to a given direction of the $n$-space, but
it could not find a consequent list of its topological properties.
Of course the shadow boundaries have been considered frequently in
convexity theory. I mention only two interesting results in context
of Baire categories see \cite{gruber-sorger} and \cite{zamfirescu}.
In \cite{gruber-sorger} the authors proved that a typical shadow
boundary under parallel illumination from a direction vector has
infinite $(n-2)$-dimensional Hausdorff measure, while having
Hausdorff dimension $(n-2)$. In \cite{zamfirescu} it is shown that,
in the sense of Baire categories, most $n$-dimensional convex bodies
have infinitely long shadow boundaries if the light vector comes
along one of $(n-2)$-dimensional subspaces.

\begin{defi}
Let $K$ be a compact convex body in $n$-dimensional Euclidean space
$E^n$ and let $S^{(n-1)}$ denote the $(n-1)$-dimensional unit sphere
in $E^n$. For ${\bf x}\in S^{(n-1)}$ the {\bf shadow boundary}
$S(K,{\bf x})$ of $K$ in direction ${\bf x}$ consists of all points
$P$ in $bd K$ such that the line $\{P+\lambda {\bf x}: \lambda \in R
\mbox{ (real numbers)} \}$ supports $K$, i.e. it meets $K$ but not
the interior of $K$. The shadow boundary $S(K,{\bf x})$ is {\bf
sharp} if any above supporting line of $K$ intersects $K$ exactly in
the point $P$. If $S(K,{\bf x})$ is not sharp, in general, it may
have {\bf sharp point} for that the above uniqueness holds.
\end{defi}

To make this paper more self-contained, we list and show some
topological properties of this set. ( Of course some of these are
well-known fact.)

\begin{statement}

The shadow boundary decomposes the boundary of $K$ into three
disjoint sets. These are $S(K,{\bf x})$ itself, moreover
\begin{eqnarray}
K^{+} &:= &\{{\bf y}\in bd K | \mbox{ there is } \tau >0 \mbox{ such
that } {\bf y}-\tau \cdot {\bf x}\in
int(K) \},\\
K^{-}&:= &\{{\bf y}\in bd K | \mbox{ there is } \tau >0 \mbox{ such
that } {\bf y}+\tau \cdot {\bf x}\in int(K) \}, \nonumber
\end{eqnarray}
\noindent respectively.
\end{statement}

We call the congruent (thus homeomorphic) sets $K^{+}$ and $K^{-}$
the {\bf positive and negative part} of $bdK$, respectively.

\proof From the definition of the sets and the convexity of $K$ it
is obvious then the above sets are disjoints, and each of the points
of $bd (K)$ belongs to one of them. On the other hand the line
through the origin and parallel to the vector ${\bf x}$ intersects
the boundary of $K$ the points $P^{+}\in K^{+}$ and $P^{-}\in K^{-}$
showing that there are not empty. \qed

\begin{defi}
We call the points $P^{+}$ and $P^{-}$ the {\bf positive} and {\bf
negative pole} of $K$, respectively. The intersection of $bd (K)$ by
a 2-plane containing the poles we call the {\bf longitudinal
parameter curve} of $K$.
\end{defi}

\begin {statement}

$S(K,{\bf x})$ is an at least $(n-2)$-dimensional connected, closed
(so compact) set in $bd(K)$, the sets $K^{+}$ and $K^{-}$ are
homeomorphic copies of $R^{(n-1)}$ giving two arc-wise connected
components of their union.

\end {statement}

\proof

Let $p_{\bf x}$ be the parallel to ${\bf x}$ projection of the
embedding space $R^n$ onto a hyperplane orthogonal to the vector
${\bf x}$. Since the orthogonal projection is a contraction then it
is continuous (i.e. a mapping of the space).  $p_{\bf x}(K)$ is a
convex body of the image hyperplane, the interior of this body is
the image of the sets $K^{+}$ and $K^{-}$, respectively and its
boundary is the image of $S(K,{\bf x})$. Since $p_{\bf x}$ on
$K^{+}$ is a bijection it is easy to see that it is a homeomorphism,
giving the first statement on $K^{+}$ and $K^{-}$, respectively. Of
course the union of their is open thus the shadow boundary is
closed.

Since $R^{(n-1)}$ is arc-wise connected the second statement on
$K^{+}$ follows from the fact that an arc connecting two points of
$K^{+}$ and $K^{-}$ decomposed into two relative open sets by
$K^{+}$ and $K^{-}$, which is impossible. (Arc-wise connectivity of
a set implies its connectivity, too.) This also means that the
shadow boundary separates the boundary of $K$ so using a theorem of
Alexandrov  (Th. 5.12 in vol.I of \cite{alexandrov}) we get that the
topological dimension of $S(K,{\bf x})$ at least $n-2$ as we stated.

We now prove that $S(K,{\bf x})$ is connected. Assume that $K_1$ and
$K_2$ two closed disjoint subsets of the shadow boundary for which
$K_1\cup K_2=S(K,{\bf x})$. First we observe that each of the metric
segments lying on a longitudinal parameter curve and parallel to
${\bf x}$ connected subset of $S(K,{\bf x})$, thus its points (by
the "basic lemma of connectivity" see vol.I p.13 in
\cite{alexandrov}) belongs either to the set $K_1$ or to the set
$K_2$. If now $C_1$ and $C_2$ the sets defined by the union of those
longitudinal parameter curves which intersect the sets $K_1$ and
$K_2$, respectively, we have the equalities $C_1\cup C_2=bd K$ and
$C_1\cap C_2=\{P^{+},P^{-}\}$, respectively. Of course the sets
$C_i$ are closed in $bd K$, meaning that the sets $C_i\setminus
\{P^{+},P^{-}\}$ gives a decomposition of $bd K \setminus
\{P^{+},P^{-}\}$ into disjoint relative closed subsets, too. Since
the latter set is connected it follows that either $K_1$ or $K_2$ is
empty. \qed

\begin {remark} In the first example of this section we
construct such a centrally symmetric convex body, which shadow
boundary is not locally connected or locally contractible. This
means that (in general) the shadow boundary does not an absolute
neighbourhood retract (ANR) meaning that it is not generalized
manifold (especially topological manifold or curved polyhedron).
\end {remark}

In general the dimension of $S(K,{\bf x})$ is $(n-2)$ or $(n-1)$. We
prove that there is a $(n-2)$-dimensional closed, connected subset
of $S(K,{\bf x})$ separating $bd K$, too.

\begin {lemma}
The boundary (frontier) of the closure of the set $K^{+}$ (denoted
by $bd(cl(K^+))$) is a closed, connected $n-2$ dimensional subset of
$S(K,{\bf x})$ separating the boundary of $K$.
\end {lemma}

\proof

By its definition it is closed. Since $ cl(K^+)\supset K^+$ and $
cl(K^+)\cap K^-=\emptyset$ we have (by Statement 1) $K^+\subset
cl(K^+)\subset cl(K^+)\cup S(K,{\bf x})$. On the other hand
$bd(cl(K^+))\cap K^+=\emptyset $ since the set $K^+$ is an open one,
then we get that $bd(cl(K^+))\subset S(K,{\bf x})$.

The separating property follows from the fact that the union of the
pairwise disjoint sets $bdK \setminus cl(K^+)$, $int (cl(K^+))$, $bd
(cl(K^+))$ fills the boundary of $K$ and the first two sets are
open. (Here $int(\cdot)$ mean the interior of the set in the
bracket.)

Now the separating property implies (again by the Alexandrov theorem
above) the inequality $dim(bd(cl(K^+)))\geq (n-2)$. On the other
hand a closed connected set of dimension $(n-1)$ on $bd K$ contains
an interior point relative to $bd K$ (see p.174 in volI. of
\cite{alexandrov} ) which contradicts to the definition of
$bd(cl(K^+))$. \qed

Before proving the main statement of this paragraph, we consider
some examples showing the possibilities of the strange attitude of
the sets, defined above.

\begin {examples}
\begin {enumerate}
\item Consider the following sequence of segments of
of $R^3$ (with respect to a fixed orthonormal coordinate system)
$s_n=\{(t,\frac{1}{n},\frac{\sqrt{n^2-1}}{n}) | -1\leq t\leq 1\},
n\in N$ with limit segment $s=\{(0,0,1) | -1\leq t\leq 1\}$. Denote
by $\cal{L}^+$ the union of these segments. If now we connect the
point $(1,\frac{1}{n},\frac{\sqrt{n^2-1}}{n})$ to the point
$(-1,\frac{1}{n+1},\frac{\sqrt{(n+1)^2-1}}{n+1})$ by the arc, which
is the intersection of the 2-plane orthogonal to the plane of the
first two coordinates through the above two points with the cylinder
of points $\{(t,r,s) | -1\leq t\leq 1, r^2+s^2=1, r,s\geq 0\}$ for
every $n$, we get a connected, closed set lies on the cylinder
$\{(t,r,s) | -1\leq t\leq 1, r^2+s^2=1\}$. This is not arc-wise nor
locally connected and it is not locally contractible one. It is easy
to see that if we add this curve to its reflected image in the
coordinate plane $[x,z]$, and the union curves we add to its
reflected image in the plane $[x,y]$, we have a centrally symmetric
closed curve belongs to the cylinder $\{(t,r,s) | -1\leq t\leq 1,
r^2+s^2=1\}$. The convex hull of this curve $\gamma$ (similar to the
so-called topologist's sine curve) is a centrally symmetric convex
body. If the direction of the light is parallel to the axe $x$, then
we have $S(K,{\bf x})=bd(cl(K^+))=bd(cl(K^-))=\gamma $. Since it is
not locally contractible it could not be an ANR. Consequently it is
not generalized (hence topological) manifold nor curved polyhedron.

\item We refer to the Example 1. in paper \cite{gho2} 
which is a presentation of a shadow boundary such like can be
seen on Fig.1. More precisely, take the parameter values
$t_{i,j}=\frac{j}{2^{i}}2\pi$, where $0\leq i$ is integer and $1\leq
j\leq 2^i$ is odd number. The diadic rational points of the circle
are the points $S_{i,j}=(\cos (t_{i,j}), sin(t_{i,j}))$ of the
subspace $E^2$ with respect to an orthonormed basis. (Note that we
define the points $S_{0,1}$ and $S_{1,1}$ in the $0^{th}$ and
$1^{th}$ steps, respectively, and -- in the $i^{th}$ step -- we
consider further $2^{i-1}$ points of form $S_{i,j}$ of the circle.)
Let now $s_{i,j}$ be a segment orthogonal to the subspace $E^2$
whose midpoint is $S_{i,j}$ and its length is equal to
$\frac{1}{2^{i-2}}$ if $i\geq 2$ and is equal to 2 if $i=0,1$. The
point sets
$$
C^*:=C\cup (\cup_{i,j}\{s_{i,j}\}) \mbox{ and } K:=conv C^*
$$
are central symmetric, here $ conv $ abreviates convex hull. This
body is also closed, see it in Fig.1. In this case $S(K,{\bf x})$ is
an ANR but is not a manifold, while the sets
$bd(cl(K^+))=bd(cl(K^-))$ are the same metric circle.

\begin{figure}
\begin{picture}(0,130)
\put(150,130){\special{em:graph fig2.pcx}}
\end{picture}
\caption{Shadow boundary which is not a topological manifold.}
\end{figure}

\item In this example $bd(cl(K^+))$ and $bd(cl(K^-))$ are the common
boundary of the sets $S(K,{\bf x})$ and $K^+$, $S(K,{\bf x})$ and
$K^-$, respectively, and homeomorphic to $S^1$ but $S(K,{\bf x})$ is
not 1-manifold nor 2-manifold with boundary. Consider the regular
octahedron as $K$ and let the direction of the light be parallel to
an edge of $K$. The shadow boundary is the polyhedron containing
four face of $K$ connecting to each other with a common edge
(parallel to {\bf x}) or a common vertex. The sharp points of the
shadow boundary are these two vertices. Of course of this points
there is no neighbourhood homeomorphic to a segment or a plane.
\end {enumerate}
\end {examples}

We concentrate now to the cases when the above sets are topological
manifolds. We can state now the followings:

\begin{theorem}
If the shadow boundary $S(K,{\bf x})$ is a topological manifold of
dimension $(n-2)$ then it is homeomorphic to the $(n-2)$-sphere
$S^{(n-2)}$. If it is an $(n-1)$-dimensional manifold with boundary
then it is homeomorphic to the cylinder $S^{(n-2)}\times [0,1]$.
($[0,1]$ means the unit interval of the real line.)
\end{theorem}

Before the proof we recall some using definitions and theorems on
the theme of cell-like mappings. There are several good papers on
this important chapter of geometric topology (e.g. \cite{lacher1},
\cite{lacher2} or \cite{semmes}). We here follow the setting up of
the nice paper of W.J.R.Mitchell and D.Repovs
\cite{mitchell-repovs}.

A non-empty compactum $K$ is said to be {\bf cell-like} if for some
embedding of $K$ in an ANR $M$, the following property holds: For
every neighbourhood $U$ of $K$ in $M$, there exists a neighbourhood
$V$ such that $K\subset V\subset U$ and the inclusion
$i:V\longrightarrow U$ is nullhomotopic. (We recall that the space
$Y$ is ANR if whenever it is embedded as a closed subset of a
separable metric space, that it is a retract of some neighbourhood
of it in that space.) Given a map (of course it is a continuous
function ) $f:X\longrightarrow Y$, we say $f$ is {\bf cell-like}, if
for each $y\in Y$, the inverse image $f^{-1}(y)$ is cell-like. We
will use the following theorem.

\begin{theorem}[Cell-like Approximation Theorem for manifolds]
Let $n\not =3$ be integer. For every cell-like map
$f:M\longrightarrow N$ between topological $n$-manifolds, and every
$\varepsilon > 0$, there is a homeomorphism $h:M\longrightarrow N$
such that $d(f,h)<\varepsilon $ in the sup-norm metric on the space
of all continuous maps (called by {\bf near homeomorphism}).
\end{theorem}

The long history of this result can be read in
\cite{mitchell-repovs}. We note that in the 3-dimensional case has
an analogous approximation theorem for a subset of the class of
cell-like mappings called by the class of ${\bf cellular}$ maps. A
set of the manifold $M$ called {\bf cellular}, if it is an
intersection of a sequence of closed cell $B_i$ of $M$ with the
properties $K\subset B_i$ and $B_{i+1}\subset B_{i}$. A map is
cellular if the inverse images are cellular sets. Cellularity
originated in the work of M.Brown \cite{brown} while the concept of
cell-likeness introduced by R.C.Lacher in \cite{lacher1}. The
concept of cellularity depends on the embedding of the examined
metric space $K$ in $M$, this dependence on embedding was eliminated
in the concept of cell-likeness. In fact, (in the manifold case)
every cellular map is a cell-like map, since every cellular set is a
cell-like one. Contrary, if we consider a wild arc in $R^3$ which
has non-simply-connected complement it is non-cellular set, while
the standard embeddings manifestly are cellular in $R^3$, showing
that it will be a cell-like set.

We also remark that a cellular (or cell-like) map in general is not
a near homeomorphism since there is a cellular map on $S^1\times
[0,1]$ to $S^1$ which are not near homeomorphism (since if we have a
near homeomorphism between compact metric space then these space are
homeomorphic to each other).

Now we can prove the main theorem above of this paragraph.

\proof [Theorem 1] Consider now again the projection $p_{\bf x}$
(was defined in the proof of Statement 2), and restrict it to the
shadow boundary of $K$ parallel to ${\bf x}$. It is of course a
cell-like map because of the point inverses are points or segments,
respectively. So for $n\not = 5$ the approximation theorem of
cell-like maps shows that this restricted map is a near
homeomorphism on $S(K,{\bf x})$ to $S^{(n-2)}$ meaning that there
are homeomorphic to each other. On the other hand this map is also
cellular, since the metric segments and points of $S(K,{\bf x})$ are
cellular sets in it. To prove this, let $s=p_{\bf x}({\bf v})$ be a
segment in $S(K,{\bf x})$  for a ${\bf v}\in S(K,{\bf x})$. If now
$Q\in s$ a point, consider a metric ball $B_{\epsilon}(Q)\subset bd
(K)$ with center $Q$ and radius $\epsilon >0$ for which
$int(B_{\epsilon}(Q))\cap S(K,{\bf x})$ homeomorphic to $R^{(n-2)}$.
Such an $\epsilon >0$ is exists. In fact, $Q$ has a neighbourhood
$N_Q$ in $S(K,{\bf x})$ homeomorphic to $R^(n-2)$. If we can choose
a point $P_\epsilon \in B_{\epsilon}(Q)\cap S(K,{\bf x})$ does not
belong to $N_Q$ for every $\epsilon$ then we have a sequence of
points $(P_\epsilon)$ having the same property and converges to $Q$.
Since $N_Q$ is open in $S(K,{\bf x})$, this is impossible. Thus
there is an $\epsilon >0$ for which $B_{\epsilon}(Q)\cap S(K,{\bf
x})=B_{\epsilon}(Q)\cap N_Q$ meaning that $int(B_{\epsilon}(Q))\cap
S(K,{\bf x})$ is an open subset of $N_Q$ relative to the topology of
$S(K,{\bf x})$. Of course $\epsilon$ depends on $Q$, but $s$ is
compact set meaning that there are finite number of points $Q_i$ and
$\epsilon _i$, for which $\cup int(B_{\epsilon ^*}(Q_i))\subset s$
(by the minimum value $\epsilon ^*$ of the $\epsilon$'s) is the
interior of the closed cell $\cup (B_{\epsilon ^*}(Q_i))\supset s$.
If we take into consideration then the property $B_{\epsilon}(Q)\cap
S(K,{\bf x})=B_{\epsilon}(Q)\cap N_Q$ also holds for every $\epsilon
'$ which less or equal to $\epsilon $, we have that the cellularity
property hold for $s$ for an infinite sequence of sets above, where
$\epsilon ^*$ tends to zero.

Observe now that if $S(K,{\bf x})$ is an $(n-1)$-manifold with
boundary (which dimension (by definition) is $(n-2)$) then its
boundary has two connected components which are equal to
$bd(cl(K^+))$) and $bd(cl(K^-))$), respectively.

First we see that $bd(cl(K^+))$) is the common boundary points of
$cl(K^+)$ and $S(K,{\bf x})$ meaning that $bd(cl(K^+))\subset
bd(S(K,{\bf x}))$). (Analogously we have that $bd(cl(K^-))\subset
bd(S(K,{\bf x}))$.

Second we note that there is no point of $int(cl(K^+))$ belonging to
$S(K,{\bf x})$ because for such a point $P$

-- either there is a neighbourhood $P\in U\subset S(K,{\bf x})$
homeomorphic the $(n-1)$-dimensional half-space implies that $P$ is
a boundary point of $cl(K^+)$ (in $U$ have points $Q$ with
neighbourhood $V\subset S(K,{\bf x})$ homeomorphic to $R^{(n-1)}$)
for which $Q\in V\subset U$ hence $Q$ is a point of the complement
of $cl(K^+)$)

-- or there is a neighbourhood $P\in U\subset S(K,{\bf x})$
homeomorphic the $n-1$-dimensional space shows that $P$ is in the
interior of $S(K,{\bf x})$ contradicting with the assumption that it
is a point of $int(cl(K^+))$. This means that $int(cl(K^+))=K^+$ and
so $bd(cl(K^+))=bd(K^+)$ is the common boundary of $K^+$ and
$S(K,{\bf x})$ hence by Lemma 1. $bd(cl(K^+))$ is a connected closed
subset of the boundary of $S(K,{\bf x})$.

Using now the fact that $bd(cl(K^-))$ is the image of $bd(cl(K^+))$
by a central projection, we have the similar properties for
$bd(cl(K^-))$, too. (It is the common boundary of $K^-$ and
$S(K,{\bf x})$.) We prove that the boundary of $S(K,{\bf x})$ is the
disjoint union of these two sets. The relation $bd(S(K,{\bf
x}))\subset bd(cl(K^-))\cup bd(cl(K^+))$ is obvious.

Consider a point $P$ from the intersection $bd(cl(K^-))\cap
bd(cl(K^+))$. Let $B$ be a metric $(n-1)$-ball around this point
with sufficiently small radius $\epsilon >0$, which intersection
with the neighbourhood $U\subset S(K,{\bf x})$ (homeomorphic to a
half-space of $R^{(n-1)}$) is the "half" of $B$ meaning that this
intersection is a topological half-space of dimension $(n-1)$ and
its complement homeomorphic to an open half-space. (We note that
such an $\epsilon >0$ and ball $B$ there exists as it can be seen
easily by similar arguments as we did it in the proof of cellularity
of a metric segments.) Since $B$ contains points from $K^+$ and
$K^-$,respectively we have a contradiction, because $K^+$ and $K^-$
is not separating by $S(K,{\bf x})$.(There is no points of $S(K,{\bf
x})$ in this complementary domain.)

This means that the boundary of $S(K,{\bf x})$ has two connected
components which are the common boundary of $S(K,{\bf x})$ and
$K^+$, $S(K,{\bf x})$ and  $K^-$, respectively. Of course these sets
are also $(n-2)$-manifolds connected with straight line segments in
all of their points. So we have $S(K,{\bf x})=bd(cl(K^+))\times
[0,1]$ holds. We have to prove only, that in this case $bd(cl(K^+))$
homeomorphic to $S^{(n-2)}$, too. Since $p_{\bf x}$ on $bd(cl(K^+))$
into $S^{(n-2)}$ is also a cell-like (and cellular) map,
$bd(cl(K^+))$ is an $(n-2)$-dimensional manifold and this restricted
map is onto one, the last statement of the Theorem follows from
$Theorem 2$, too. \qed

\section{General parameter spheres}

We now recall the definition of general parameter spheres. (see
\cite{gho2} ).
\begin{defi}
Let
$$
\lambda _0:=\sup\{t | tK\cap (tK+{\bf x})=\emptyset \}
$$ be the smallest value $\lambda $ for which $\lambda K$ and $\lambda K+{\bf x}$ are
intersect. Then the {\bf generalized parameter sphere of $K$}
corresponding to the direction {\bf x} and to the parameter $\lambda
\geq \lambda _0$ is the following set:
$$
\gamma _{\lambda }(K,{\bf x}):=\frac{1}{\lambda }(bd(\lambda K)\cap
bd(\lambda (K)+{\bf x})).
$$
\end{defi}

In \cite{gho2} we mentioned that in general the above sets are not
topological spheres of dimension $(n-2)$ and are not homeomorphic to
each other. For example the dimension of $\gamma _{\lambda
_0}(K,{\bf x})$ may be $0$, $1$ $\cdots$ $(n-1)$ while the
topological dimension of $\gamma _{\lambda }(K,{\bf x})$ is at least
$(n-2)$ because it divides the surface of $K$. We remark that the
interiors of the given two caps of the boundary are also
homeomorphic to each other as in the case of shadow boundary. In
fact a centrally projection from $\frac{1}{\lambda }{\bf x}$ sending
the left half of $bd K$ onto the left one of
$\frac{1}{\lambda}(bd(\lambda K)+{\bf x})$ is an appropriate
homeomorphism. (This latter set is congruent to the right half of
$bd K$ by the central symmetry of $\lambda K\cap \lambda K+{\bf
x}$.) Also we proved that the shadow boundary $S(K,{\bf x})$ is the
limit of the generalized parameter spheres $\gamma _{\lambda
}(K,{\bf x})$, with respect to the Haussdorff metric, when $\lambda
$ tends to infinity.

We also saw (in the proof of Lemma 1 in \cite{gho2}) that the
general parameter sphere $\gamma _{\lambda }(K,{\bf x})$ is the
shadow boundary of the convex body $\frac{1}{\lambda}(\lambda K\cap
\lambda K+{\bf x})$ thus the statements of the previous section can
be adapted to their. The purpose of this section to examine the
manifold case, we prove the following two statement:

\begin{theorem}
I, The shadow boundary $S(K,{\bf x})$ is an $(n-2)$-dimensional
manifold if all of the non-degenerated parameter spheres $\gamma
_{\lambda }(K,{\bf x})$ with $\lambda >\lambda _0$ are
$(n-2)$-dimensional manifolds, contrary if $S(K,{\bf x})$ is an
$(n-2)$-dimensional manifold then all of the general parameter
spheres are ANR.

II, The shadow boundary $S(K,{\bf x})$ $(n-1)$-dimensional manifold
with boundary iff there is a $\lambda$ for which the general
parameter sphere $\gamma _{\lambda }(K,{\bf x})$ is an
$(n-1)$-dimensional manifolds with boundary.
\end{theorem}

Before the proof we recall a nice theorem of M.Brown on the
projective limit of compact metric space and corresponding near
homeomorphisms (see \cite{brown 2} or \cite{uspenskij}). A map from
$X$ to $Y$ between compact metric spaces is a {\bf near
homeomorphism} if it is in the closure of the set of all
homeomorphisms from $X$ onto $Y$, with respect to the sup-norm
metric on the space $C(X,Y)$ of all maps from $X$ to $Y$. Now the
mentioned theorem is:

\begin{theorem}[M.Brown]
Let $(X_n)$ be an inverse sequence of compact metric spaces with
limit $X_{\infty }$. If all bonding maps $X_k\longrightarrow X_n$
are near homeomorphisms, then so are the limit projections
$X_k\longrightarrow X_{\infty }$.
\end{theorem}

Let us give an example showing that we need distinguish the above
two cases of the theorem.

\begin{example}
Consider the union of the connecting rectangles $\pm \{(r,1,t) |
-1\geq r,t\geq 1\}$, $\pm \{(r,s,t) | r+s=2, 0\geq r\geq 2, -1\geq
t\geq 1\}$, $\pm \{(r,s,t) | r-s=2, 0\geq r\geq 2, -1\geq t\geq 1\}$
and the segments $\pm \{(r,0,2) | -\frac{3}{2}\geq r\geq
\frac{3}{2}\}$. The convex hull $K$ of this set is a convex
polyhedron. If now the vector ${\bf x}$ is the position vector
directed into the point $(4,0,0)$ we have three important values for
the parameters of the generalized parameter spheres. For $\lambda
_0=1$ the degenerated sphere $\gamma _{\lambda _0}(K,{\bf x})$ is a
segment. For $1<\lambda \leq \frac{5}{4}$ the general parameter
spheres $\gamma _{\lambda }(K,{\bf x})$ are homeomorphic to $S^1$.
On the range $\frac{5}{4}<\lambda \leq \frac{3}{2}$ the general
parameter sphere $\gamma _{\lambda }(K,{\bf x})$ is a simplicial
complex containing one or two-dimensional simplices, respectively.
(This space is an ANR but is not a topological manifolds.) Finally,
in the last parameter domain $\lambda > \frac{3}{2}$ the $\gamma
_{\lambda }(K,{\bf x})$ are homeomorphic to the cylinder $S^1\times
[0,1]$. Since $S(K,{\bf x})$ is the union of six quadrangles,
parallel to the axe $x$ it is also a cylinder.
\end{example}

We also remark that if $S(K,{\bf x})$ is an $(n-2)$-dimensional
manifold than probably all of the non-degenerated parameter spheres
are also. Unfortunately we can not prove this statement.

\proof[Theorem 3] First we note that -- for every $\lambda
_0<\lambda'<\infty $ -- $S(K,{\bf x})$ can be got as the inverse
limit space $X_{\infty }$ of the metric spaces $X_\lambda:=\gamma
_{\lambda }(K,{\bf x})$ for $\lambda'<\lambda $. In fact, by Lemma 1
in\cite{gho2} for $\lambda
>\lambda_0$ the intersection of $\gamma _{\lambda }(K,{\bf x})$ by a
longitudinal parameter curve, say $r$ is a segment then $r\cap
\gamma _{\mu }(K,{\bf x})$ with $\mu>\lambda$ is also a segment
containing the segment $r\cap \gamma _{\lambda }(K,{\bf x})$. So in
this case the union of the set $r\cap \gamma _{\mu }(K,{\bf x})$ is
the segment $r\cap S(K,{\bf x})$. On the other hand we have two
possibilities in the case when $r\cap \gamma _{\lambda }(K,{\bf x})$
is a point. In the first one $r\cap S(K,{\bf x})$ is a point, too,
meaning that for all $\mu
>\lambda$ for $r\cap \gamma _{\lambda }(K,{\bf x})$ also is a point.
If now $r\cap S(K,{\bf x})$ is a segment then we have value $\lambda
' >\lambda $ with the property that if $\mu >\lambda '$ then $r\cap
\gamma _{\mu }(K,{\bf x})$ is a segment, too. In this latter case
$r\cap S(K,{\bf x})=\cup\{r\cap \gamma _{\mu }(K,{\bf x})\}$. Define
now the left end of a segment parallel to ${\bf x}$ as that end of
it which has less parameter in the usual parametrization with
respect to ${\bf x}$ (Meaning that the point of a line parallel to
${\bf x}$ is in the form $P+\tau {\bf x}$ with a point $P$ of its
line.) The bounding map $p_{\lambda, \mu}$ for $\gamma _{\mu
}(K,{\bf x})$ to $\gamma _{\lambda }(K,{\bf x})$ ($\mu
>\lambda $) let us define in the following way:

For a point $P$ of $\gamma _{\mu }(K,{\bf x})$
$$
p_{{\lambda, \mu}}(P)=\left\{
\begin{array}{ll}
r\cap \gamma _{\lambda }(K,{\bf x}) & \mbox{ if $r\cap \gamma
_{\lambda }(K,{\bf x})$ is a point}\\

P & \mbox{ if $r\cap \gamma _{\lambda }(K,{\bf x})$ is a segment and
$P\in r\cap \gamma _{\lambda }(K,{\bf x})$}\\

\mbox{the left end of $r\cap \gamma _{\lambda }(K,{\bf x})$} &
\mbox{ if $P\in r\cap \gamma _{\mu }(K,{\bf x})\setminus r\cap
\gamma _{\lambda }(K,{\bf x})$}
\end{array}
\right. .
$$

The continuity of this function (with respect to theirs relative
metric) is obvious and the inverse (projective) limit space
$X_{\infty}$ can be identified by $S(K,{\bf x})$ with the limit
mappings $p_{\mu}$ defined an analogously way from $S(K,{\bf x})$ to
$\gamma _{\mu }(K,{\bf x})$ as the above functions $p_{\lambda,
\mu}(P)$. (Of course, we have the sufficient equality $p_{\mu_p',\mu
''}\circ p_{\mu '}=p_{\mu ''}$ for $\mu ''>\mu '$.)

Using the theorems 2 and 4 above, the proof of the first direction
of the first statement is follows easily. In fact, if for $\lambda
>\lambda _0$ the space $\gamma _{\lambda }(K,{\bf x})$ is an
$(n-2)$-manifold then using Theorem 2 we know that the bounding maps
$p_{\mu', \mu''}: \gamma _{\mu '' }(K,{\bf x})\longrightarrow \gamma
_{\mu' }(K,{\bf x})$ are near homeomorphisms. By Theorem 4 now we
get that the limit projections $p_{\lambda}$ are also near
homeomorphisms. This means that the space $S(K,{\bf x})$ is also an
$(n-2)$ manifold.

Conversely, if now $S(K,{\bf x})$ is an $(n-2)$-dimensional manifold
than it is locally contractible. By Lemma 1 in \cite{gho2} this
also means that all of the general parameter spheres are locally
contractible, too. On the other hand the general parameter spheres
can be considered as the compact subsets of $E^{(n-1)}$ meaning that
there are ANR. (See Theorem 8 in p.117 in \cite{daverman}).

The proof of both of the second statement uses now Theorem 1. If
first we have a general parameter sphere $\gamma _{\lambda }(K,{\bf
x})$ which is $(n-1)$-dimensional manifold with boundary, then by
Theorem 1 it is a cylinder with boundaries homeomorphic to
$S^{(n-2)}$. In this case the shadow boundary contains this general
parameter spheres showing that every point-inverses with respect to
$p_{\bf x}$ are segments (with non-zero lengthes). On the other
hand, the sets $bdK^+\cap S(K,{\bf x})$ and $bd K^+\cap \gamma
_{\lambda }(K,{\bf x})$ are agree, showing that $S(K,{\bf x})$ is a
cylinder based on an $(n-2)$ manifold homeomorphic to $S^{(n-2)}$.
Since $bd K^-\cap S(K,{\bf x})$ homeomorphic to $S^{(n-2)}$ (by
central symmetry) and these two sets are disjoints we have that in
fact, $S(K,{\bf x})$ homeomorphic to $S^{(n-2)}\times [0,1]$ as we
stated.

Conversely, if $S(K,{\bf x})$ is an $(n-1)$-manifold with boundary,
then it is (by Theorem 1) homeomorphic to $S^{(n-2)}\times [0,1]$.
Since this cylinder is compact there is a positive value
$\varepsilon$ less or equal to the length of each of the segments
intersected from the shadow boundary by a longitudinal parameter
curve. This means that there is a $\lambda <\infty $ such that
$\gamma _{\lambda }(K,{\bf x})\subset S(K,{\bf x})$. The
intersection $\gamma _{\lambda }(K,{\bf x})\cap K^+$ is the same as
the intersection $S(K,{\bf x})\cap K^+$ which is one of the two
components of the boundary of $S(K,{\bf x})$ homeomorphic with
$S^{(n-2)}$. For this $\lambda $ it is possible a trivial
point-inverse with respect the map $p_{\bf x}$ as we saw in the
example of this section, but for every $\lambda'>\lambda $ the
general parameter sphere $\gamma _{\lambda' }(K,{\bf x})$ is a
cylinder. Using now the fact that it is also the shadow boundary of
a central symmetric convex body which positive part is the set
$K^+$, we have that it is also manifold with boundary homeomorphic
to $S^{(n-2)}\times [0,1]$.

\section{On the bisector and its embedding}

In this section we investigate the bisector $H_{\bf x}$ (or
equidistant set of the starting and ending point of the vector {\bf
x}) using the system $\lambda \gamma _{\lambda }(K,{\bf x})$ of
compact metric spaces. Our goal proving the following theorem:

\begin{theorem}
$H_{\bf x}$ is an $(n-1)$-dimensional manifold if and only if the
non-degenerated general parameter spheres $\gamma _{\lambda }(K,{\bf
x})$ are manifolds of dimension $(n-2)$.
\end{theorem}

Since the neighbourhoods (with respect to $H_{\bf x}$ of the point
$\frac{1}{2}{\bf x}$ can not be homeomorphic to either $R^n$ or a
half space, this is the only manifold case for $H_{\bf x}$.

\proof To prove the first direction we use Theorem 1. From this we
know that the general parameter spheres are homeomorphic copies of
$S^{(n-2)}$, respectively. Construct now the bisector $H_{\bf x}$ as
the disjoint union of the sets $\gamma _{\lambda }(K,{\bf x})$ for
$\lambda \geq \lambda_0$. The set $H_{{\bf x},\mu}=\{ \lambda \gamma
_{\lambda }(K,{\bf x}) | \mu \geq \lambda\geq \lambda_0\}$ obviously
homeomorphic to $\gamma _{\lambda }(K,{\bf x})\cup K_{\lambda}^+$
meaning that it is a homeomorphic copy of the closed
$(n-1)$-dimensional ball. Thus $int H_{{\bf x},\mu}$ is homeomorphic
to $R^{n-1}$ for each $\mu\geq \lambda_0$. Applying now a M.Brown
theorem (see in \cite{uspenskij} or \cite{brown}) saying that if a
topological space is the union an increasing sequence of open
subsets, each of which is homeomorphic to $R^{(n-1)}$, then it is
also homeomorphic to $R^{(n-1)}$, we get the required result.

Conversely, the projection $p_{\bf x}:H_{\bf x}\longrightarrow
R^{(n-1)}$ is a cellular map between two manifolds with the same
dimension if $H_{\bf x}$ is homeomorphic to $R^{(n-1)}$. Thus it is
a near homeomorphism meaning that the restriction of it to the
compact metric space $\lambda \gamma _{\lambda }(K,{\bf x})$ is a
near homeomorphism, too. But its image is the boundary of a convex
compact $(n-1)$-dimensional body gives at once that it is a
homeomorphic copy of $S^{(n-2)}$. Hence the general parameter
spheres $\gamma _{\lambda }(K,{\bf x})$ for $\lambda >\lambda_0$ are
manifold of dimension $(n-2)$ as we stated. \qed

\corollary The proof of the first direction of the conjecture follows from 
the theorems 1,3 and 5. In fact, if $H_{\bf x}$ is a topological hyperplane then each of the  non-degenerated general parameter spheres are a homeomorphic copy of $S^{(n-2)}$ by Theorem 5 and Theorem 1. So by Theorem 3 we get that the shadow boundary is also a homeomorphic copy of $S^{(n-2)}$ which is the statement of the mentioned direction of our conjecture. 

On the other hand we could prove in Theorem 3 only that if $S(K,{\bf x})$ is a homeomorphic copy of $S^{(n-2)}$ then the non-degenerated parameter spheres are ANR, respectively, thus the manifold property for the bisector does not follows immediately from our theorems. Furthermore, in the manifoldian case we prove only that the bisector is a homeomorphic copy of $R^{(n-1)}$ so a weaker property as the required one. So we have to investigate the question of embedding. In fact, all of the example in geometric topology giving a non-standard (wild) embedding of
a set into $R^n$ based on the observation that the
connectivity properties of the complement (with respect to $R^n$) of
the set can change applying a homeomorphism to it. In our case, the
complement for example of the bisector (which is now the
homeomorphic copy of $R^{(n-1)}$) is the disjoint union of the
homeomorphic copies of $R^n$ giving the chance the existence of a
homeomorphism on $R^n$ to itself sending the bisector to a
hyperplane. It is a well-known fact that a manifold homeomorphic to $S^{(n-1)}$ in $S^n$
is unknotted if and only the closure of the components of its complement are
homeomorphic copy of the closed n-cells $B^n$. This means that in the manifold case the embedding of the
shadow boundary and the general parameter spheres are always standard, respectively meaning the existence of a homeomorphism of the boundary of $K$ into itself sending these sets into a standard $(n-1)$-dimensional sphere of $bd K$. In the case of the bisector we have to prove a little bit more precisely. Let $\varphi$ be a homeomorphism sending $H_{\bf x}$ into $R^{(n-1)}$ (which is now a hyperplane $H$ of $R^n$). We consider the compactification of the embedding space by an element denoted by $\infty$. Extend first the map $\varphi$ to the  compact space $H_{\bf x}\cup \{\infty\}$ by the condition $\varphi (\infty)=\infty $. Of course this extended map gives a homeomorphism between the sets $H_{\bf x}\cup \{\infty\}$ and $R^{(n-1)}\cup \{\infty\}$. Since the closure of the components of the complement of $H_{\bf x}\cup \{\infty\}$ in $R^n\cup \{\infty\}$ are closed n-cells the homeomorphism $\varphi $ can be extended to a homeomorphism $\Phi :R^n\cup\{\infty\}\longrightarrow R^n\cup\{\infty\}$. Since by our method we have: $\Phi (\infty)=\varphi(\infty)=\infty $ and $\Phi (H_{\bf x})=H$ we get that the bisector is a topological hyperplane as we stated. Thus the following statement is fulfill:
\begin{theorem}
In the manifold case the embeddings of $H_{\bf x}$, $S(K,{\bf x})$ and $\gamma_{\lambda}(K,{\bf x})$ are standard, respectively, meaning that if the bisector homeomorphic to $R^{(n-1)}$ then it is a topological hyperplane.
\end{theorem}

\end{document}